\documentclass{amsart}

\usepackage{amsmath}
\usepackage{amsfonts}
\usepackage{amssymb}
\usepackage{amsthm}
\usepackage{graphics}
\usepackage [dvips]{epsfig}
\usepackage{amssymb,amsfonts,amstext,amsmath,graphicx,epic,amsthm,color,times}

\usepackage{graphicx}
\usepackage{t1enc}

\newtheorem{theorem}{Theorem}[section]
\newtheorem{lemma}[theorem]{Lemma}
\newtheorem{corollary}[theorem]{Corollary}
\newtheorem{proposition}[theorem]{Proposition}
\theoremstyle{definition}
\newtheorem{definition}[theorem]{Definition}

\theoremstyle{remark}
\newtheorem{remark}[theorem]{Remark}

\numberwithin{equation}{section}

\usepackage[all]{xy}
\title[Flux geometry and dynamical systems]{On symplectic dynamics}

\author{St\'ephane Tchuiaga}
\address{Department of Mathematics of the University of Buea, 
 South West Region, Cameroon}
\email{tchuiagas@gmail.com}

\begin{document}

\renewcommand\contentsname{Table of contents}
\renewcommand\refname{References}
\renewcommand\abstractname{Abstract}
\pagestyle{myheadings}
\markboth{MICHAEL MAGEN}{ESPACE DES ORBITES D'UNE GRASSMANNIENNE}

\maketitle

\begin{abstract}

This paper continues to carry out a foundational study of Banyaga's 
topologies of a closed symplectic manifold \cite{Ban08a}. 
Our intension in writing this paper is to provide 
several symplectic analogues of some  results found in the study of Hamiltonian dynamics. 
Especially, 
without appealing to the positivity of the symplectic displacement energy, we point out the impact of 
 the $L^\infty$ version of Banyaga's Hofer-like metric in the investigation of the symplectic nature of
 the $C^0-$limit of a sequence of symplectic maps. This result is the symplectic analogue of 
a result that was proved by Hofer-Zehnder \cite{Hof-Zen94}
 (for compactly supported Hamiltonian diffeomorphisms on
 $\mathbb{R}^{2n}$), and then reformulated by Oh-M\"{u}ller \cite{Oh-M07} for Hamiltonian diffeomorphisms in general. 
Furthermore, we extend to symplectic isotopies the regularization procedure for Hamiltonian 
paths introduced by Polterovich \cite{Polt93}, and then we use it to 
prove the equality between the two versions of Banyaga's Hofer-like norms defined on the identity 
component in the group of symplectomorphisms. This result was announced in \cite{BanTchu}. It shows the uniqueness of Banyaga's 
Hofer-like geometry, and then yields the symplectic analogue of a 
result that was proved by Polterovich \cite{Polt93}. Finally, we elaborate the symplectic analogues 
of some approximation results found in Oh-M\"{u}ller \cite{Oh-M07}, and make some remarks on flux theory.\\


\end{abstract}

\section{Introduction}\label{SC00}

The Hofer geometry started with the remarkable paper of Hofer \cite{Hofer90} that 
introduced the Hofer topologies on the space of Hamiltonian diffeomorphisms of a symplectic 
manifold. 
These topologies motivated various studies in the field of Hamiltonian  
dynamics. In particular, Hofer-Zehnder \cite{Hof-Zen94} derived almost all the basic 
formulae and some perspectives for the subsequent development of Hamiltonian dynamics.
 A thorough discussion of Hofer topologies can be found in \cite{Hof-Zen94, Lal-McD95, Oh-M07,Polt93}.\\ 
Recently, it was shown by Banyaga \cite{Ban08a} that the Hofer topologies admit 
a natural generalization to the set of all time-one maps of symplectic isotopies of a symplectic manifold. 
In particular, if a symplectic manifold  is such that the identity component in its group of symplectic diffeomorphisms is reduced to the group 
of all Hamiltonian diffeomorphisms, then 
 the corresponding Banyaga's topologies \cite{Ban08a} reduce to Hofer's topologies. \\
These facts attest that it is judicious to investigate whether the analogues of 
some results found in the study of Hamiltonian dynamics can be elaborated 
in the context of Banyaga's Hofer-like geometry or not.\\ This motivated the results of the present
paper:  
In Section 2, we recall some fundamental facts concerning symplectic mappings and isotopies. Section 2.8 
introduces Hopf-Rinow theorem from Riemannian geometry and shows its implication in the study of Hofer's norms 
with respect to a certain class of functions. In Section 2.9, using Hodge's theory, we 
show that  
Polterovich's regularization process for Hamiltonian isotopies admits a natural generalization to
 symplectic isotopies.  Section 3 deals 
with the main results of the present paper. 
Here, studying Banyaga's topologies, we use the results of Sections 2.8 and 2.9 to 
prove that Banyaga's Hofer like-geometry is invariant 
under the choice of  Banyaga's Hofer-like norm. We show an impact of the $L^\infty$ version of the Hofer-like 
metric in the investigation of the symplectic nature of a homeomorphism which is the $C^0-$limit of a sequence of symplectic diffeomorphisms. 
This 
follows by combining Hodge's decomposition theorem 
of symplectic isotopies together with the standard continuity theorem of ODE for Lipschitz vector fields. Furthermore, we 
prove that if a loop is homotopic (relatively to a fixed base point) to 
a closed Hamiltonian orbit, then the symplectic area swept by the latter under the symplectic flow generated by any non exact closed $1-$form 
is trivial. This can be viewed (in a certain sense) as the dual form   
of a following well-known result from flux geometry. Section 4 contains some approximation lemmas that generalize  
some results found in the study of Hamiltonian dynamics.

\section{Preliminaries}\label{SC0}
Let $M$ be a smooth closed manifold of dimension $2n$. A differential $2-$form $\omega$ on $M$ is called a symplectic form if
 $\omega$ is closed and nondegenerate. In particular, any symplectic manifold is oriented.  
From now on, we shall always assume that $M$ admits a symplectic form $\omega$. A diffeomorphism $\phi:M\rightarrow M$ 
is called symplectic if it preserves the symplectic form $\omega$.
\subsection{Symplectic vector fields}
The symplectic structure $\omega$ on $M,$ being nondegenerate, induces an isomorphism between vector fields $Z$ and $1-$forms on $M$ given by 
$Z\mapsto \omega(Z,.)=: \iota(Z)\omega$. 
A vector field $Z$ on $M$ is symplectic if $\iota(Z)\omega$ is closed. In particular, 
a symplectic vector field $Z$ on $M$ is said to be  a Hamiltonian vector 
field if $\iota(Z)\omega$ is exact. It follows from the definition of symplectic vector fields that, 
 if the first de Rham cohomology group of the manifold $M$ is trivial (i.e. $H^1(M,\mathbb{R}) = 0$), 
then all the symplectic vector fields induced  
by a symplectic form $\omega$ on $M$ are Hamiltonian. If we equip $M$ with a Riemannian metric $g$ 
(any differentiable manifold $M$ can be equipped with a Riemannian metric), then any harmonic $1-$form $\alpha$ on $M$ 
determines a symplectic
 vector field $Z$ such that $\iota(Z)\omega = \alpha$ (so-called harmonic vector field, see \cite{Ban08a}). 
In view of Hodge's theory,  a sufficient condition that guarantees the existence of nontrivial 
harmonic vector fields on a symplectic manifold $(M, \omega)$ is that $H^1(M,\mathbb{R}) \neq 0$. 
Note that $H^1(M,\mathbb{R})$ is a topological invariant, i.e. 
it 
 does not depend on the differentiable structure on $M$ 
  and depends 
 only on the underlying topological structure of $M$ \cite{FW71}.
\subsection{Symplectic isotopies}
An isotopy $\{\phi_t\}$  of the symplectic manifold $(M,\omega)$ is said to be symplectic if for each $t$, the vector field 
$Z_t = \dfrac{d\phi_t}{dt}\circ\phi_t^{-1}$ is symplectic. In particular 
 a symplectic isotopy $ \{\psi_t\}$ is called Hamiltonian if for each $t$,  
$Z_t = \dfrac{d\psi_t}{dt}\circ\psi_t^{-1}$ is Hamiltonian, i.e. there exists a smooth 
function $F: [0,1]\times M \rightarrow \mathbb{R}$ called Hamiltonian such that
 $\iota(Z_t)\omega = dF_t$. As we can see, any Hamiltonian isotopy determines a 
 Hamiltonian $F: [0,1]\times M \rightarrow \mathbb{R}$ up to an additive constant . 
Throughout the paper we assume that all Hamiltonians are normalized in the following way: given a Hamiltonian $F: [0,1]\times M \rightarrow \mathbb{R}$
we require that $\int_M F_t \omega^n = 0$. We denote by $\mathcal{N}([0,1]\times M\ ,\mathbb{R})$
the space of all smooth normalized Hamiltonians and by $Ham(M,\omega)$ 
the set of all time-one maps of Hamiltonian isotopies. If we equip $M$ with a Riemannian metric 
$g$, then a symplectic isotopy $ \{\theta_t\}$ is said to be harmonic if for each $t$,
 $Z_t = \dfrac{d\theta_t}{dt}\circ\theta_t^{-1}$ is harmonic. We denote 
by $Iso(M, \omega)$ the group of all symplectic isotopies of $(M,\omega)$ and by $Symp_0(M,\omega)$ 
the set of all time-one maps of symplectic isotopies.

\subsection{Harmonics $1-$forms}
 From now on, we assume that $M$ is equipped with a Riemannian metric $g$, 
and denote by $\mathcal{H}^1(M,g)$ the space of harmonic $1-$forms on $M$ 
with respect to the Riemannian metric $g$. In view of the Hodge theory,  $\mathcal{H}^1(M,g)$ 
is a finite dimensional vector space over $\mathbb{R}$ which is isomorphic to $H^1(M,\mathbb{R})$ (see \cite{FW71}). 
The dimension of $\mathcal{H}^1(M,g)$ is the first Betti number of the manifold $M$, denoted by 
$b_1$. Taking ${(h_i)_{1\leq i\leq b_1}}$ as a basis of the vector space $\mathcal{H}^1(M,g)$, we 
 equip $\mathcal{H}^1(M,g)$ with the Euclidean norm $|.|$ defined as follows : 
for all $H$ in $\mathcal{H}^1(M,g)$ with $H = \Sigma_{i = 1}^{b_1}\lambda_ih_i$ 
we have $$|H| :=\Sigma_{i = 1}^{b_1}|\lambda_i|.$$ It is convenient to compare the above Euclidean norm with the well-known 
uniform sup norm of differential $1-$forms. For this purpose, let's recall the definition of the uniform sup norm of a differential $1-$form
 $\alpha$ on $M$. For all $x\in M$, we know that $\alpha$ induces a linear map $\alpha_x : T_xM \rightarrow \mathbb{R}$ 
whose norm is given by $$\|\alpha_x\| = \sup\{ |\alpha_x(X)| : X\in T_xM,  \|X\|_g = 1\},$$
where $\|.\|_g$ is the norm induced 
on each tangent space $ T_xM$ by the Riemannian metric $g$. Therefore, the uniform sup norm of $\alpha$, say $|.|_{0}$ is defined by 
 $|\alpha|_{0} = \sup_{x\in M}\|\alpha_x\|.$
 In particular, when $\alpha$ is a harmonic $1-$form (i.e $\alpha = \Sigma_{i = 1}^{b_1}\lambda_ih_i$), we obtain the following estimates, 
$$|\alpha|_{0} \leq \Sigma_{i = 1}^{b_1}|\lambda_i| | h_i|_0\leq E|\alpha|,$$ 
where $ E := \max_{1\leq i\leq b_1}|h_i|_{0}$. 
If the basis ${(h_i)_{1\leq i\leq b_1}}$ is such that $E \textgreater 1$, then one can always normalize such a basis so that  
 $E$ equals $1$. Otherwise, the identity $|\alpha|_{0}\leq E|\alpha|$ reduces to $|\alpha|_{0}\leq|\alpha|$.
We denote by $ \mathcal P\mathcal{H}^1(M,g)$, the space of smooth mappings $\mathcal{H}: [0,1]\rightarrow\mathcal{H}^1(M,g)$.

\subsection{A description of symplectic isotopies \cite{BanTchu}}\label{SC3}
In this subsection, from the group of symplectic isotopies, we shall deduce another group which
will be convenient later on (see \cite{BanTchu}). Consider $\{\phi_t\}$ to be a symplectic isotopy, for each $t$, the vector field 
$Z_t = \dfrac{d\phi_t}{dt}\circ(\phi_t)^{-1}$ satisfies $d\iota(Z_t)\omega =0$. So, it follows from   
 Hodge's theory that
$\iota(Z_t)\omega$ decomposes as the sum of an exact $1-$form $dU_t^\Phi$ and a harmonic $1-$form $\mathcal{H}_t^\Phi$ (see \cite{FW71}). 
Denote by $U$ the Hamiltonian $U^\Phi = (U^\Phi_t)$ normalized,  
and by $ \mathcal{H}$ the smooth family of harmonic $1-$forms $ \mathcal{H}^\Phi = (\mathcal{H}^\Phi_t)$. 
In \cite{BanTchu}, the authors denoted by $\mathfrak{T}(M, \omega, g)$ the Cartesian product 
$\mathcal{N}([0,1]\times M ,\mathbb{R})\times \mathcal P\mathcal{H}^1(M,g)$, and equipped it with a
 group structure which makes the bijection
\begin{equation}
  Iso(M, \omega)\rightarrow\mathfrak{T}(M, \omega, g), \Phi\mapsto (U, \mathcal{H})
\end{equation}
 a group isomorphism. Denoting the map just constructed by $\mathfrak{A}$, the authors denoted  
any symplectic isotopy $\{\phi_t\}$ as $\phi_{(U,\mathcal{H})}$ 
to mean that the mapping $\mathfrak{A}$ maps $\{\phi_t\}$ onto $(U,\mathcal{H})$, and $(U,\mathcal{H})$ is 
called the 
``generator'' of the symplectic path $\phi_{(U,\mathcal{H})}$. In particular, any symplectic isotopy of the form 
$\phi_{(0,\mathcal{H})}$ is considered to be a harmonic isotopy, while any 
symplectic isotopy of the form $\phi_{(U,0)}$ is considered to be a Hamiltonian 
isotopy. The product in $\mathfrak{T}(M, \omega, g)$ is given by, 
\begin{equation}
  (U,\mathcal{H})\Join(V ,\mathcal{K}) = ( U + V\circ\phi_{(U,\mathcal{H})}^{-1} 
+ \widetilde \Delta(\mathcal{K},\phi_{(U,\mathcal{H})}^{-1}), \mathcal{H} + \mathcal{K} )
\end{equation}
The inverse of $(U,\mathcal{H})$, denoted $\overline{(U,\mathcal{H})}$ is given by 
\begin{equation}
 \overline{(U,\mathcal{H})} = (- U\circ\phi_{(U,\mathcal{H})} - \widetilde\Delta(\mathcal{H},\phi_{(U,\mathcal{H})}),-\mathcal{H})
\end{equation}
where for each $t$, $\phi_{(U,\mathcal{H})}^{-t} :=(\phi_{(U,\mathcal{H})}^{t})^{-1}$, and 
$\widetilde \Delta_t(\mathcal{K},\phi_{(U,\mathcal{H})}^{-1})$ is the function\\ 
$\Delta_t(\mathcal{K},\phi_{(U,\mathcal{H})}^{-1}) 
:= \int_0^t\mathcal{K}_t(\dot{\phi}_{(U,\mathcal{H})}^{-s})\circ\phi_{(U,\mathcal{H})}^{-s}ds$ normalized.\\ 

Here is
 a consequence of the Hodge decomposition theorem of symplectic isotopies \cite{Ban08a}. 

\begin{proposition}\label{mainremark} Every $(U,\mathcal{H})\in\mathfrak{T}(M, \omega, g)$ decomposes in a unique way as  
\begin{equation}
 (U,\mathcal{H}) = (0,\mathcal{H})\Join(U\circ\phi_{(0,\mathcal{H})}, 0)
\end{equation}
\end{proposition}
{\it Proof}.  
 Let $\{\phi_t\}$ be the symplectic isotopy generated by $(U,\mathcal{H})$. It follows from \cite{Ban08a} that
 $\{\phi_t\}$ decomposes in 
a unique way as $\{\phi_t\} =  \{\rho_t\}\circ \{\psi_t\}$ where $\{\rho_t\}$ is a harmonic isotopy  
and $\{\psi_t\}$ is a Hamiltonian isotopy. Now, for each $t$, we compute $\dot\phi_t = 
\dot \rho_t + (\rho_t)_\ast(\dot\psi_t)$ 
and derive that 
$$\iota(\dot\phi_t)\omega = \iota(\dot\rho_t)\omega + (\rho_t^{-1})^\ast(\iota(\dot\psi_t)\omega) = \mathcal{H}_t + dU_t.$$
It follows immediately from the above identities that $\iota(\dot\rho_t)\omega = \mathcal{H}_t$ and 
$ \iota(\dot\psi_t)\omega = (\rho_t)^\ast(dU_t)$ for each $t$, i.e.    
$ \{\rho_t\}$ is generated by $(0,\mathcal{H})$ and $ \{\psi_t\}$ is generated by $(U\circ \{\rho_t\}, 0)$ so that    
$(U,\mathcal{H}) 
= (0,\mathcal{H}) \Join (U\circ \{\rho_t\}, 0).$ The uniqueness of this decomposition is supported by 
the uniqueness of the Hodge decomposition of symplectic isotopies \cite{Ban08a}. This completes the proof. $ \hspace{11.4cm}\Box$\\

\subsection{Reparameterization of symplectic isotopies \cite{BanTchu, Oh-M07}}\label{RPSC4}
We shall need the following basic formula for the  generator of a reparameterized
symplectic path. Let $\Phi = \{ \phi_t\}\in Iso(M,\omega)$ which is generated by $( U,\mathcal{H})$, and let   
$\xi :[0,1]\rightarrow[0,1] $ be a smooth function. The reparameterized path 
$\Phi^{\xi}: t\mapsto \phi_{\xi(t)}$ is generated by the
 element $( U,\mathcal{H})^{\xi}$ defined by  
\begin{equation}\label{Rep}
 ( U,\mathcal{H})^{\xi} = ( U^{\xi},\mathcal{H}^{\xi}),
\end{equation}
 where $\mathcal{H}^{\xi} $ is the smooth map $t\mapsto\dot\xi(t)\mathcal{H}_{\xi(t)}$, 
while $ U^{\xi}$ is the smooth map 
$t\mapsto  \dot\xi(t)U_{\xi(t)}$, and 
$\dot\xi(t)$ is the derivative of $\xi$ with respect to $t$. 

\begin{definition}(\cite{Oh-M07}). 
 Given a curve $\xi :[0,1]\rightarrow \mathbb{R}$, its norm $\|\xi\|_{ham}$ is defined by 
$\|\xi\|_{ham} = \|\xi\|_{C^0} + \|\dot\xi\|_{L_1}$ where 
$\|\dot\xi\|_{L_1} = \int_0^1|\dot\xi(t)|dt,$ and \hspace{0.1cm}$\|\xi\|_{C^0} = \sup_t|\xi(t)|.$
\end{definition}
\subsection{Boundary flat symplectic isotopies \cite{BanTchu, Oh-M07}}\label{SC6}
\begin{definition}(\cite{BanTchu}). 
  Let $( U,\mathcal{H})\in\mathfrak{T}(M, \omega, g)$. We say that $( U,\mathcal{H})$ is boundary flat if there exists $\delta\in]0,1[$ such that 
$( U_t,\mathcal{H}_t) = (0,0)$ for all $t\in[0,\delta[\cup]1-\delta,1].$
\end{definition}

It follows from the definition above that a path $ \{ \phi_t\}$ is boundary flat  
if there exists a constant $0\textless\delta\textless1$ such that 
 $\phi_t = id$ for all $0\leq t\textless\delta$ and $\phi_t = \phi_1$ for all $1-\delta\textless t\leq 1$.\\

The following results show some properties of the Hamiltonians $\Delta(\mathcal{H},\Phi)$.

\begin{proposition}\label{SC41}
 Let $(M,\omega)$ be a closed  symplectic manifold. If $\Phi = \{\phi_t\}\in Iso(M,\omega)$ is Hamiltonian, then 
for any $\mathcal{H}\in \mathcal P\mathcal{H}^1(M,g)$, 
the Hamiltonian $\Delta(\mathcal{H},\Phi)$ is normalized, i.e. $\int_M \Delta_t(\mathcal{H},\Phi)\omega^n = 0$ for each $t$.
\end{proposition}

{\it Proof}. Assume that $\Phi = \{\phi_t\}$ is Hamiltonian and for each $t$, set  $Z_t = \dfrac{d\phi_t}{dt}\circ\phi_t^{-1}$. 
We use the expression $ \Delta_t(\mathcal{H},\Phi) = \int_0^t \mathcal{H}_t(Z_s)\circ\phi_sds$ to obtain 
\begin{equation}\label{PHPP}
 \int_M \Delta_t(\mathcal{H},\Phi)\omega^n  
 = \int_0^t(\int_M \mathcal{H}_t(Z_s)\omega^n)ds,
\end{equation}
for each $t$, but the identity $\mathcal{H}_t\wedge \omega^n = 0$ yields 
$(\iota(Z_{s})\mathcal{H}_t)\wedge\omega^n + 
(\iota(Z_{s})\omega^n)\wedge \mathcal{H}_t = 0$ for all $s\in[0,t]$, and this implies in turn that  
\begin{equation}\label{PHP}
 \int_M \mathcal{H}_t(Z_{s})\omega^n = - \int_M (\iota(Z_{s})\omega^n)\wedge \mathcal{H}_t
= -n\int_Md[U_s\circ\phi^{-1}_s\omega^{n-1}\wedge\mathcal{H}_t],
\end{equation}
 for all $s\in[0,t]$, where $U$ denotes the generating Hamiltonian of the path $ \Phi$. Combining (\ref{PHPP}) 
and (\ref{PHP}) yield, 
\begin{equation}\label{PPHP}
 \int_M \Delta_t(\mathcal{H}, \Phi)\omega^n = \int_0^t(\int_M \mathcal{H}_t(Z_s)\omega^n)ds 
= -n\int_Md[(\int_0^t U_s\circ\phi^{-1}_sds)\omega^{n-1}\wedge\mathcal{H}_t],
\end{equation}
Applying Stokes' theorem in the right-hand side of (\ref{PPHP}) leads to 
\begin{equation}
 \int_M \Delta_t(\mathcal{H},\Phi)\omega^n = -n\int_{\partial M} (\int_0^t U_s\circ\phi^{-1}_sds)\omega^{n-1}\wedge\mathcal{H}_t =0
\end{equation}
 since $\partial M = \emptyset$. This completes the proof. $ \hspace{6.5cm}\Box$

\begin{lemma}\label{cP11}
Let $\Phi_1 = (\phi_1^t)$ and $\Phi_2 = (\phi_2^t)$ be two symplectic isotopies. Let $\mathcal{H}$ be a smooth family of closed $1-$forms. Then,  
for each $t$, there exists a constant $C$ which depends on $\Phi_1$, $\Phi_2$, and $t$ such that, 
$$\Delta_t(\mathcal{H},\Phi_1 \circ\Phi_2 ) =  \Delta_t(\mathcal{H},\Phi_2) + \Delta_t(\mathcal{H},\Phi_1)\circ\phi_2^t + C.$$
\end{lemma}

{\it Proof.} For a fixed $t$, we have, 
$$
\begin{array}{lllllllccccccc}
 d\Delta_t(\mathcal{H}, \Phi_1 \circ\Phi_2) &=& (\phi_1^t \circ\phi_2^t)^\ast\mathcal{H}_t - \mathcal{H}_t\cr\cr
&=& (\phi_2^t)^\ast((\phi_1^t)^\ast(\mathcal{H}_t)) - \mathcal{H}_t\cr\cr
&=& (\phi_2^t)^\ast(\mathcal{H}_t + d\Delta_t(\mathcal{H},\Phi_1)) - \mathcal{H}_t \cr\cr
 &= & (\phi_2^t)^\ast(\mathcal{H}_t) + d \Delta_t(\mathcal{H},\Phi_1)\circ\phi_2^t - \mathcal{H}_t\cr\cr
&=& \mathcal{H}_t +  d\Delta_t(\mathcal{H},\Phi_2) + d \Delta_t(\mathcal{H},\Phi_1)\circ\phi_2^t - \mathcal{H}_t\cr\cr
 &=& d\Delta_t(\mathcal{H},\Phi_2) + d \Delta_t(\mathcal{H},\Phi_1)\circ\phi_2^t.
\end{array}
$$
 It follows from the above 
estimates that : $$\Delta_t(\mathcal{H},\Phi_1 \circ\Phi_2 )  = 
\Delta_t(\mathcal{H},\Phi_2) +  \Delta_t(\mathcal{H},\Phi_1)\circ\Phi_2 + C.$$ This achieves the proof. $\square$ 

\subsection{The $C^0-$metric}\label{SC2}
 Let $Homeo(M)$ be the homeomorphisms' group of $M$ equipped with the $C^0-$ compact-open topology. This is the 
metric topology induced by the distance
  $$d_0(f,h) = \max(d_{C^0}(f,h),d_{C^0}(f^{-1},h^{-1})),$$
where $d_{C^0}(f,h) =\sup_{x\in M}d (h(x),f(x))$ and 
 $ d$ is a distance on $M$ induced by the Riemannian metric $g$. 
 On the space of all continuous paths $\varrho:[0,1]\rightarrow Homeo(M)$ such that $\varrho(0) = id$, 
we consider the $C^0-$topology as the metric topology induced by the metric  
$$\bar{d}(\lambda,\mu) = \max_{t\in [0,1]}d_0(\lambda(t),\mu(t)).$$ 

\subsection{Scholium}\label{HOPFRINOW}
Here, using Hopf-Rinow theorem, we show some bounded properties of the Hamiltonians $ \Delta(\mathcal{H},\Phi)$ with respect 
to the Hofer-norms. This will be useful later on.\\

When the manifold $M$ admits a complete Riemannian metric, 
by Hopf-Rinow theorem one can choose the path $\gamma$ to be a geodesic, and
 its length is bounded from above by the 
diameter of the manifold $M$, i.e. $\int_0^1\|\dot{\gamma}(s)\|_gds\leq diam(M)$ 
where $diam(M)$ denotes the diameter of $M$ with respect to the Riemannian metric $g$.
 Let $\mathcal{H} = (\mathcal{H}_t)\in  \mathcal P\mathcal{H}^1(M,g)$
 and $\{\phi_t\}$ be an isotopy. 
Since for each $t$, the harmonic $1-$form $ \mathcal{H}_t$ is closed, it follows that 
\begin{equation}\label{Sholum}
 \phi_t^\ast(\mathcal{H}_t) - \mathcal{H}_t  = d \Delta_t(\mathcal{H},\{\phi_t\}),
\end{equation}
for each $t$. It follows from (\ref{Sholum}) that 
\begin{equation}\label{Mos}
 u_t(x):=\int_{\gamma}\phi_t^\ast( \mathcal{H}_t) - \mathcal{H}_t =\Delta_t( \mathcal{H}, \{\phi_t\})(x) -
\Delta_t(\mathcal{H},\{\phi_t\})(x_0),
\end{equation}
for each $t$. 
In view of the above fact, to study the behavior of the Hamiltonian 
$(x,t)\mapsto\Delta_t(\mathcal{H} ,\{\phi_t\})(x)$ with respect to the Hofer topologies we will only need  
 to study 
 the behavior of its associated Hamiltonian $(x,t)\mapsto u_t(x)$ with respect to the uniform sup norm since  
the two norms are equivalent. For instance, 
let $y_0$ be any point of $M$ that realizes the supremum of the function $x\mapsto|u_t(x)|$. 
We derive from the triangle inequality that, 
$$\sup_x|u_t(x)|\leq |\int_{\gamma_{y_0}} \mathcal{H}_t| + |\int_{\gamma_{y_0}}\phi_t^\ast(\mathcal{H}_t)|
$$ $$\leq 
|\mathcal{H}_t|\int_0^1\|\dot{\gamma}_{y_0}(s)\|_g ds +   \sup_{t,s}|D\phi_{t}( \gamma_{y_0}(s))||\mathcal{H}_t|\int_0^1\|\dot{\gamma}_{y_0}(s)\|_gds,$$
$$\leq diam(M)(1 + \sup_{t,s}|D\phi_{t}( \gamma_{y_0}(s))|)| \mathcal{H}_t|,$$
where $D\phi_t$ is the tangent map of $\phi_t$. This yields the following estimates 
\begin{equation}\label{Mos1}
 \int_0^1osc(\Delta_t(\mathcal{H} ,\{\phi_t\}))dt 
\leq 2diam(M)(1 + \sup_{t,s}|D\phi_{t}( \gamma_{y_0}(s))|) \int_0^1|\mathcal{H}_t|dt.
 \end{equation}
\begin{equation}\label{Mos2}
 \max_{t\in [0,1]} osc(\Delta_t(\mathcal{H} ,\{\phi_t\}))
\leq 2 diam(M)(1 + \sup_{t,s}|D\phi_{t}( \gamma_{y_0}(s))|) \max_{t\in [0,1]}|\mathcal{H}_t|.
 \end{equation}

\subsection{Regularization of symplectic isotopies}\label{SC7} 
A regularization process for Hamiltonian paths is due to Polterovich \cite{Polt93}.
As far as we know a general regularization process for symplectic isotopies is unknown. The pathern is underlying to the proof of Lemma 3.3 found 
in \cite{BanTchu}. Here, we use Hodge's theory to show that Polterovich's regularization process for Hamiltonian isotopies admits a natural generalization
 to symplectic isotopies. 
Before we start, note that a symplectic path $\{\phi_t\}$
is said to be regular if for every $t$, the tangent vector $\dot\phi_t$ to the path $\{\phi_t\}$ does not vanish. For  instance, 
let $\Phi = \phi_{( U,\mathcal{H})}\in Iso(M,\omega)$, in view of a Polterovich's result (\cite{Polt93},Proposition $5.2.A$), for
 the above Hamiltonian $U$, there exists a Hamiltonian loop $\phi_{( r,0)}$ which is
 close to the constant loop identity (in the $C^\infty-$sense), and in particular, its generating function $r$ 
is arbitrarily small in the $L^{(1,\infty)}$ version of Hofer's norm so that
\begin{equation}\label{Polterovich}
 osc(-r_t + U_t) \neq0
\end{equation}
 for all $t$. Now, consider 
 $( V,\mathcal{K})$ to be the product $\overline{(r,0)}\Join(U,\mathcal{H}),$ which can be written immediately as  
\begin{equation}\label{tPolt}
 ( V,\mathcal{K}) = 
(-r\circ\phi_{( r,0)} + U\circ\phi_{( r,0)} + \widetilde\Delta(\mathcal{H},\phi_{( r,0)}), \mathcal{H}).
\end{equation}
We claim that the isotopy generated by the element $( V, \mathcal{K})$ just constructed above is regular. As a matter of fact, 
 assume 
that there exists a time $s$ for which the vector field $ X_s = \dot{\phi}_{( V, \mathcal{K})}^s$ vanishes identically, and this is 
equivalent to,
\begin{equation}\label{tPol}
 \iota(X_s)\omega = dV_s + \mathcal{K}_s = 0.
\end{equation}
Inserting (\ref{tPolt}) in (\ref{tPol}), we obtain
\begin{equation}\label{ReG}
 d(-r_s\circ\phi_{( r,0)}^s + U_s\circ\phi_{( r,0)}^s + \Delta_s(\mathcal{H},\phi_{( r,0)})) + \mathcal{H}_s = 0.
\end{equation}
From (\ref{ReG}), it follows that 
the harmonic $1-$form $\mathcal{H}_s$ is exact, and then the latter must be trivial because in view of 
Hodge's theory \cite{FW71}, any exact harmonic form is trivial. This suggests that
$\Delta_s(\mathcal{H},\phi_{( r,0)}) = 0$, which implies in turn that the function   
$x\mapsto (-r_s\circ\phi_{( r,0)}^s + U_s\circ\phi_{( r,0)}^s)(x)$ is constant. This contradicts the assumption 
$ osc(r_t -U_t) \neq 0$ for all $t$, and the claim follows. 

As a consequence of the above regularization process, we derive that using 
any regular symplectic path $\phi_{( V, \mathcal{K})}$ 
we can define a function $\zeta :[0,1]\rightarrow[0,1]$ to be the 
inverse of the map,  $$ s\mapsto \dfrac{\int_0^s(osc(V_t) + |\mathcal{K}_t|  )dt}{\int_0^1( osc(V_t) + |\mathcal{K}_t| )dt},$$
and the derivative of $\zeta$ is given explicitly by :
\begin{equation}\label{FNB}
  \zeta'(s) = \frac{\int_0^1( osc(V_t) + |\mathcal{K}_t|) dt}{osc(V_{\zeta(s)}) + |\mathcal{K}_{\zeta(s)}|}.
\end{equation}
Note that if $\zeta$ is only $C^1$,   
then we can approximate $\zeta$ in the $C^1-$topology by a smooth diffeomorphisms  $\kappa:[0,1]\rightarrow[0,1]$ 
that fixes $0$ and $1$ (see \cite{Hirs76}).
 This will enable us to prove the uniqueness of Banyaga's Hofer-like metric.

\section{Main results}
Throughout this section, we introduce the main results of this paper.\\

We shall need the definitions of 
Banyaga's  topologies \cite{Ban08a}.\\

According to \cite{Ban08a}, 
the so-called $L^{(1,\infty)}$ version and $L^{\infty}$ version of Banyaga's Hofer-like lengths of any  
$\Phi = \phi_{(U,\mathcal{H})}\in Iso(M,\omega)$ are defined respectively by, 
\begin{equation}\label{blg1}
 l^{(1,\infty)}(\Phi) = \int_0^1 osc(U_t) + |\mathcal{H}_t|dt, 
\end{equation}
\begin{equation}\label{blg2}
 l^\infty(\Phi) 
= \max_{t\in[0,1]}( osc(U_t) + |\mathcal{H}_t|).
\end{equation}
Clearly $l^{(1,\infty)}(\Phi) \neq l^{(1,\infty)}(\Phi^{-1})$ unless $\Phi$ is Hamiltonian. Indeed,  $\Phi = \phi_{(U ,\mathcal{H})}$ 
implies that  $\Phi^{-1}
= \phi_{\overline{(U,\mathcal{H})}}$ where $\overline{(U,\mathcal{H})} = 
(-U\circ\Phi - \widetilde\Delta( \mathcal{H},\Phi),-\mathcal{H})$. 
Hence, we see that the mean oscillation of the function $U$ can be different from that of the function
$-U\circ\Phi - \widetilde\Delta( \mathcal{H},\Phi)$. But, if $\Phi$ is Hamiltonian, i.e. 
$\Phi = \phi_{(U , 0)}$, then the mean oscillation of $ U $ is equal to that
of $-U\circ\Phi$, i.e. $l^{(1,\infty)}(\Phi) = l^{(1,\infty)}(\Phi^{-1})$.  
Similarly, we have $l^{\infty}(\Phi) \neq l^{\infty}(\Phi^{-1})$ unless $\Phi$ is Hamiltonian.\\  

Before we start, let's revisit an interesting result from Hamiltonian dynamics.
 In the special case that $Ham(M,\omega) = Symp_0(M,\omega)$ (or  $H^1(M,\mathbb{R}) = 0$), 
it is known that the Hofer length for Hamiltonian paths   
 have an impact in the investigation of 
the Hamiltonian nature of a homeomorphism which is the $C^0-$limit of a sequence of Hamiltonian diffeomorphisms  (\cite{Hof-Zen94}, Theorem 6).\\
 When $H^1(M,\mathbb{R}) \neq 0$, analyzing the proof given by Hofer-Zehnder \cite{Hof-Zen94}, it turns out
 that in the presence of the positivity result of the symplectic displacement energy, one 
can prove a symplectic analogue of Theorem 6 found by Hofer-Zehnder \cite{Hof-Zen94} 
(in the context of Banyaga's Hofer-like geometry). 
This shall follow closely to the proof given by Hofer-Zehnder \cite{Hof-Zen94}. The above arguments lead to the following delicate question.
In the lack of the positivity result of the symplectic displacement energy, does it make sense 
to think of the symplectic analogue of Theorem 6 found by Hofer-Zehnder \cite{Hof-Zen94} 
in the context of Banyaga's Hofer-like geometry (when $H^1(M,\mathbb{R}) \neq 0$)?\\

One main theorem of this paper gives an affirmative answer to the above question.
More precisely,  without appealing to the  
positivity of a symplectic displacement energy, we prove the following main theorem.

\begin{theorem}\label{maint}
 Let $(M,\omega)$ be a closed  
symplectic manifold.
 Let $\Phi_i = \{\phi_i^t\}$ be a sequence of symplectic isotopies,
 $\Psi = \{\psi^t\}\in Iso(M,\omega)$, and $\phi : M\rightarrow M$ be a map, such that
\begin{itemize}
 \item  $(\phi_i^1)$ converges 
uniformly to $\phi$, and 
\item $l^\infty(\Phi_i^{-1}\circ\Psi)\rightarrow0,i\rightarrow\infty$.
\end{itemize}
 Then we must have $\phi = \psi^1.$
\end{theorem}
This result shows not only an advantage of the $L^\infty$ Hofer-like length over 
the $L^{(1,\infty)}$ Hofer-like length but it yields the symplectic analogue 
of Theorem 6 found by Hofer-Zehnder \cite{Hof-Zen94} (in the $L^\infty$ context). In fact, the choice of
 the $L^\infty$ Hofer-like length is supported by the following 
results.

\begin{lemma}\label{maintt}  Let $\rho_i$ be a sequence of harmonic isotopies generated by $(0,\mathcal{H}^i)$ and 
let $\rho $ be another harmonic isotopy generated by $(0,\mathcal{H})$ such that 
$\max_{t\in [0,1]}|\mathcal{H}^i_t -\mathcal{H}_t |\rightarrow0,i\rightarrow\infty.$ Then the following properties  hold
\begin{enumerate}
 \item $l^\infty(\rho_i^{-1}\circ\rho)\rightarrow0,i\rightarrow\infty,$
\item $\rho_i$ converges in $\bar d$ to $\rho$.
\end{enumerate} 
\end{lemma}

 {\it Proof}. For (2), we define a sequence $(Z_t^i)$ of smooth family of harmonic vector fields by setting 
$\iota(Z^i_t)\omega = \mathcal{H}^i_t$ for each $i$ and for all $t$. Similarly, we define a smooth family 
$(Z_t)$ of harmonic vector fields by setting  
$\iota(Z_t)\omega = \mathcal{H}_t$ 
 for all $t$. Since by assumption we have $$\max_{t\in [0,1]}|\mathcal{H}^i_t -\mathcal{H}_t |\rightarrow0,i\rightarrow\infty,$$
it turns out that the sequence  
$(Z_t^i)$ converges uniformly to $(Z_t)$. Therefore,  it follows from the standard
 continuity theorem of ODE for Lipschitz vector fields that the sequence of paths generated by $(Z_t^i)$ 
must converge uniformly to the path generated by $(Z_t^i)$, i.e. $\rho_i$
 converges uniformly to $\rho$. For (1), we compute 
$$\overline{(0,\mathcal{H}^i)}\Join (0,\mathcal{H}) = (\widetilde\Delta( \mathcal{H} - \mathcal{H}^i,\rho_i),  \mathcal{H} - \mathcal{H}^i), $$
for each $i$, and we derive that to complete the proof, we only need to prove that 
$$\max_{t\in [0,1]} osc(\widetilde\Delta_t( \mathcal{H} - \mathcal{H}^i,\rho_i))\rightarrow0,i\rightarrow\infty.$$ For this purpose, 
we use (\ref{Mos2}) to derive that 
 $$\max_{t\in [0,1]} osc(\Delta_t( \mathcal{H} - \mathcal{H}^i,\rho_i))
\leq 2 diam(M)(1 + \sup_{t,s}|D\rho_i^t( \gamma_{y_0}(s))|) \max_{t\in [0,1]}|\mathcal{H}_t - \mathcal{H}^i_t|,$$
where $D\rho_i^t$ stands for the tangent map of $\rho_i^t$
for each $i$, $y_0\in M$ and $\gamma_{y_0}$ is a geodesic such that $\gamma_{y_0}(1) = y_0$ (see Section 2.8 of the present paper). 
The right-hand side in the above estimate tends to zero when $i$ goes to the infinity since the quantity 
 $(1 + \sup_{t,s}|D\rho_i^t( \gamma_{y_0}(s))|)$ is bounded for each $i$, and 
by assumption we have $\max_{t\in [0,1]}|\mathcal{H}^i_t -\mathcal{H}_t |\rightarrow0,i\rightarrow\infty.$
 This completes the proof. $\hspace{9.2cm}\Box$

\begin{remark}\label{RMTR} From Lemma \ref{maintt}, it is clear that 
if $\rho_i$ is a sequence of harmonic isotopies generated by $(0,\mathcal{H}^i)$ and 
 $\rho $ another harmonic isotopy generated by $(0,\mathcal{H})$, then 
the convergence $l^\infty(\rho_i^{-1}\circ\rho)\rightarrow0,i\rightarrow\infty$ 
is equivalent to the convergence $\max_{t\in [0,1]}|\mathcal{H}^i_t -\mathcal{H}_t |\rightarrow0,i\rightarrow\infty$.
\end{remark}

The following result shows another interesting advantage of the $L^\infty$ Hofer-like metric. 

 \begin{corollary}\label{cortmaintt} Let $\Phi_i$ be a sequence of symplectic isotopies and let $\Psi$ be a symplectic isotopy  
 such that $l^\infty(\Phi_i^{-1}\circ\Psi)\rightarrow0,i\rightarrow\infty$. 
If $\mu_i$ is the sequence of Hamiltonian paths arising 
in the Hodge decomposition of $\Phi_i$, and $\mu$ the Hamiltonian path arising 
in the Hodge decomposition of $\Psi$, then $l^\infty(\mu_i^{-1}\circ\mu)\rightarrow0,i\rightarrow\infty.$
\end{corollary}

{\it Proof}. Assume that for each $i$,  $\Phi_i$ is 
generated by $(U^i,\mathcal{H}^i)$ and  
$\Psi$ generated by $(U,\mathcal{H})$. In view of Proposition \ref{mainremark}, we have to prove that 
$\max_tosc (U^i_t\circ\phi_{(0,\mathcal{H}^i)}^{t} - U_t\circ\phi_{(0,\mathcal{H})}^t )\rightarrow0,i\rightarrow\infty.$
 For each $i$, we compute 
$$osc (U^i_t\circ\phi_{(0,\mathcal{H}^i)}^{t} - U_t\circ\phi_{(0,\mathcal{H})}^t )\leq 
osc (U^i_t - U_t )
+ osc (U_t\circ\phi_{(0,\mathcal{H}^i)}^t - U_t\circ\phi_{(0,\mathcal{H})}^t ),$$
 for all $t$, but by assumption we have $\max_t(osc (U^i_t - U_t))\rightarrow 0, i\rightarrow\infty$, while 
the uniform continuity of the map $ (t,x)\mapsto U_t(x)$ together with Lemma \ref{maintt} yield  
$$ \max_t(osc (U_t\circ\phi_{(0,\mathcal{H}^i)}^t - U_t\circ\phi_{(0,\mathcal{H})}^t ))\rightarrow 0, i\rightarrow\infty.$$ 
This completes the proof.$ \hspace{9cm}\Box$\\

{\it Proof of Theorem \ref{maint}.} Since we use not the positivity result 
of the symplectic displacement energy but the Hodge decomposition theorem of symplectic isotopies, 
the proof of this main result is rather delicate. We shall proceed in three steps. 
\begin{itemize}
 \item Step (1).{\bf(Convergence of symplectic isotopies)}. 
For each $i$, let $\rho_i = \{\rho_i^t\}$ and $\rho = \{\rho^t\}$ denoting respectively 
the harmonic isotopies arising in the Hodge decompositions of the paths $ \Phi_i = \{\phi_i^t\}$ and $\Psi = \{\psi^t\}$. 
Since by assumption, we have $l^\infty(\Phi_i^{-1}\circ\Psi)\rightarrow0,i\rightarrow\infty$, it turns 
out that $l^\infty(\rho_i^{-1}\circ\rho)\rightarrow0,i\rightarrow\infty$. This together with 
Lemma \ref{maintt}-(2) (or Remark \ref{RMTR}) tell us that 
 the sequence $\rho_i$ converges in $\bar d$ to $\rho$.

\item Step (2).{\bf(Decomposition of the map $\phi = \lim_{C^0}(\phi_i^1)$)}. 
For each $i$, let $\mu_i = \{\mu^t_i\}$ and $\mu = \{\mu^t\}$ denoting respectively 
the Hamiltonian isotopies arising in the Hodge decompositions of the paths $\Phi_i = \{\phi_i^t\}$ and $\Psi = \{\psi^t\}$. 
By assumption, the sequence
 of time-one maps $\phi_i^1$ converges uniformly to $\phi$, and in view of step $(1)$ the sequence
of time-one maps $\rho_i^1$ converges uniformly to the time-one map $\rho^1$. The preceding arguments suggest that 
the sequence of time-one maps $\mu_i^1$ converges 
uniformly to a  homeomorphism $\sigma$ since $\mu_i^1 = (\rho_i^1)^{-1}\circ\phi_i^1 $ for each $i$. 
 For instance, we compute 
$$d_{C^0}(\phi, \rho^1\circ\sigma)\leq  d_{C^0}(\phi, \rho_i^1\circ\mu_i^1) + d_{C^0}(\rho_i^1\circ\mu_i^1, \rho_i^1\circ\sigma ) 
+ d_{C^0}(\rho_i^1\circ\sigma, \rho^1\circ\sigma ),$$
 for each $i$,
and derive  that $\phi = \rho^1\circ\sigma$ since by assumption, $  d_{C^0}(\phi, \rho_i^1\circ\mu_i^1) 
= d_{C^0}(\phi, \phi_i^1)\rightarrow0,i\rightarrow\infty,$ 
and from the bi-invariance of the metric $d_{C^0}$ it follows that 
$$d_{C^0}(\rho_i^1\circ\mu_i^1, \rho_i^1\circ\sigma ) = d_{C^0}(\mu_i^1,\sigma )\rightarrow0,i\rightarrow\infty,$$ 
$$d_{C^0}(\rho_i^1\circ\sigma, \rho^1\circ\sigma ) =  d_{C^0}(\rho_i^1, \rho^1)\rightarrow0,i\rightarrow\infty.$$ 
\item Step (3).{\bf(The Hamiltonian nature of the map $\sigma$)}. 
 To achieve the proof, 
all we have to show is that $\sigma = \mu^1$ where  $\mu^1$ is the time-one map of the Hamiltonian path $\mu = \{\mu^t\}$.  
 Arguing indirectly, we find that there exists a small non-empty closed ball $B\subset M$ which is completely displaced by  
$\sigma^{-1}\circ\mu^1$, i.e. $B\cap[\sigma^{-1}\circ\mu^1](B) = \emptyset.$ Since $B$ is compact and the convergence 
$\mu_i^1\rightarrow \sigma$ is uniform, we must have $B\cap[(\mu_i^1)^{-1}\circ\mu^1](B) = \emptyset$ 
 for all sufficiently large $i$. The above arguments tell us that we can apply the 
energy-capacity inequality theorem from \cite{Lal-McD95} to obtain
\begin{equation}\label{EEQ}
 0\textless C(B)/2\leq l^{\infty}(\mu_i^{-1}\circ\mu),
\end{equation}
for all sufficiently large $i$, where $C(B)$ represents the Gromov area of the ball $B$. 
But in view of Corollary \ref{cortmaintt}, the right-hand side in (\ref{EEQ}) tends to zero when $i$ goes to the infinity. 
This contradicts the assumption $0\textless C(B)/2$. Therefore,  $\sigma = \mu^1$, and 
this yields $\phi = \rho^1\circ\mu^1 = \psi^1$. This completes the proof.$ \hspace{8.5cm}\Box$
\end{itemize}

The following result is an immediate consequence of Theorem \ref{maint}. It can justify the definition 
of strong symplectic isotopies in the $L^\infty$ context \cite{BanTchu, Ban10c, Tchuia03}.
 \begin{corollary}\label{Comaint}
 Let $(M,\omega)$ be a closed  
symplectic manifold.
Let $\Phi_i = \{\phi_i^t\}$ be a sequence of symplectic isotopies,
 let $\Psi = \{\psi_t\}\in Iso(M,\omega)$, and $\eta:t\mapsto \eta_t$ a family of maps $\eta_t : M\rightarrow M$, 
such that $\Phi_i$ converges  in $\bar d$ 
to $\eta$ and  $l^\infty(\Phi_i^{-1}\circ\Psi)\rightarrow0,i\rightarrow\infty$. Then $\eta = \Psi.$
\end{corollary}
{\it Proof.} Assume the contrary that  $\Psi \neq \eta$, i.e. 
there exists $t\in ]0,1]$ such that $\eta_t\neq\psi_t.$ Then the sequence 
of symplectic paths $\varXi_i:s\mapsto \phi_i^{st}$ contradicts Theorem \ref{maint}. This completes the proof.$ \hspace{10cm}\Box$ 

\subsection{Banyaga's Hofer-like norms}
Before we continue with further investigation of Banyaga topologies, we shall need the following notions.\\

Let $\phi\in Symp_0(M,\omega)$, using the above Banyaga's lengths introduced above,
Banyaga \cite{Ban08a} defined respectively
the  $L^{(1,\infty)}$ energy and $L^{\infty}$ energy of $\phi$ by,
\begin{equation}\label{bener1}
 e_0(\phi) = \inf(l^{(1,\infty)}(\Phi)),
\end{equation}
\begin{equation}\label{bener2}
 e^{\infty}_0(\phi) = \inf(l^\infty(\Phi)),
\end{equation}
where the infimum are taken over all symplectic isotopies $\Phi$ with time-one map equal to $\phi$. 
Therefore,  
the $L^{(1,\infty)}$ Banyaga's Hofer-like norm and the $L^{\infty}$ Banyaga's Hofer-like norm of $\phi$ are respectively defined by,
\begin{equation}\label{bny1}
 \|\phi\|_{HL}^{(1,\infty)} = (e_0(\phi) + e_0(\phi^{-1}))/2,
\end{equation}
\begin{equation}\label{bny2}
 \|\phi\|_{HL}^\infty = (e^\infty _0(\phi) + e^\infty _0(\phi^{-1}))/2.
\end{equation}
Each of the norms $ \|.\|_{HL}^\infty $ and $\|.\|_{HL}^{(1,\infty)}$ 
generalizes the Hofer norms for Hamiltonian diffeomorphisms in the following sense : 
in the special case of a closed symplectic manifold $(M,\omega)$ for which $Ham(M,\omega) = Symp_0(M,\omega)$ (or  $H^1(M,\mathbb{R}) = 0$),
 the norm $ \|.\|_{HL}^\infty $ reduces to a norm $ \|.\|_{H}^\infty $ called the $L^\infty$ Hofer norm, while 
the norm $\|.\|_{HL}^{(1,\infty)}$ reduces to a norm $ \|.\|_{H}^{(1,\infty)} $ called the $L^{(1,\infty)}$ Hofer
 norm. But, a result that was proved by Polterovich \cite{Polt93} shows that the above
 Hofer's norms are equal in general, i.e. $ \|.\|_{H}^{(1,\infty)}=  \|.\|_{H}^{\infty} $. In other words 
the norms $ \|.\|_{HL}^\infty $ and $\|.\|_{HL}^{(1,\infty)}$ are equal when $Ham(M,\omega) = Symp_0(M,\omega)$ (or  $H^1(M,\mathbb{R}) = 0$).
However, when $  Symp_0(M,\omega) \backslash  Ham(M,\omega) \neq \emptyset$ (or  $H^1(M,\mathbb{R}) \neq0 $), it is unknown whether 
the norms $ \|.\|_{HL}^\infty $ and $\|.\|_{HL}^{(1,\infty)}$ are equal or not. This motivated 
 the following main lemma. 
\begin{lemma}\label{ch2l1} Let $(M, \omega)$ be a closed symplectic manifold. For every 
 $\phi\in Symp_0(M,\omega)$, we have 
$$ \|\phi\|_{HL}^\infty = \|\phi\|_{HL}^{(1,\infty)}.$$
 
\end{lemma}

 This result was announced in Banyaga-Tchuiaga \cite{BanTchu} without any explicit proof. It 
yields the symplectic analogue of a result which is due to Polterovich (\cite{Polt93}, Lemma 5.1.C). Its  
proof is based on the following lemma which is a refined version of Lemma 3.3 found in \cite{BanTchu}. 

\begin{lemma}\label{4ch2l4} Let $(M,\omega)$ be a closed  
symplectic manifold. 
 Let  $\Phi$ be a symplectic isotopy, and let $\epsilon$ be a positive real number. Then, there exists  
 $\Psi$ be a symplectic isotopy with the same extremities than $\Phi$ which is regular such that $l^\infty(\Psi) \textless l^{(1,\infty)}(\Phi) + \epsilon.$
\end{lemma}
{\it Proof of Lemma \ref{ch2l1}}. The inequality $\|.\|_{HL}^{(1,\infty)}\leq \|.\|_{HL}^{\infty}$ 
 is clear from the definition of the energies. For the converse,  
let $ \phi\in Symp_0(M,\omega)$, we derive from the characterization of the infimum 
that for all positive real number $\epsilon$, 
we can find a symplectic path 
$\Phi_\epsilon$ that connects $\phi $ to the identity such that
$
  l^{(1,\infty)}(\Phi_\epsilon)  \leq e_0(\phi) + \epsilon.
$
But, Lemma \ref{4ch2l4} shows that there exists $\Psi_\epsilon\in Iso(M,\omega)$ with the 
same extremities than $\Phi_\epsilon$ such that 
$l^\infty(\Psi_\epsilon) \textless l^{(1,\infty)}(\Phi_\epsilon) + \epsilon.$ This yields
 $e_0^\infty(\phi) \leq l^\infty( \Psi_\epsilon) 
\textless e_0(\phi) + 2\epsilon,$ i.e.  $ e_0^\infty(\phi) \textless e_0(\phi) + 2\epsilon. $  
Similarly, we use Lemma \ref{4ch2l4} to derive that  $
 e_0^\infty(\phi^{-1})\textless e_0(\phi^{-1}) + 2\epsilon.$
Therefore, we summarize the above estimates to get
\begin{equation}\label{equality}
  \|\phi\|_{HL}^\infty = (e_0^\infty(\phi^{-1}) + e_0^\infty(\phi))/2 
 \textless (e_0(\phi) + e_0(\phi^{-1}) + 4\epsilon)/2  
\leq \|\phi\|_{HL}^{(1,\infty)} + 2\epsilon.
\end{equation}
Since (\ref{equality}) holds for all arbitrary positive $\epsilon$, we conclude that 
 $ \|\phi\|_{HL}^\infty \leq  \|\phi\|_{HL}^{(1,\infty)}$. This completes the proof.$ \hspace{9.5cm}\Box$\\

\subsection{Proof of Lemma \ref{4ch2l4}} 
In the rest of this paper, we will always denote by $r(g)$, the injectivity radius of a Riemannian metric $g$.
We will need the following result.

\begin{lemma}\label{S1C42}Let $(M,g)$ be a closed oriented Riemannian manifold. Let $\mathcal{H}\in \mathcal P\mathcal{H}^1(M,g)$.
The following facts hold : 
\begin{enumerate}
 \item Let $\Psi = \{\psi_t\}$ be an isotopy, 
and let $\xi_1,\xi_2 :[0,1]\rightarrow[0,1] $ be two smooth monotone functions that fix $0$. Then 
there exists a constant  $B_2$ which depends on $\mathcal{H}$ and $\Psi$ such that, 
$$\int_0^1osc(\Delta_t(\mathcal{H},\Psi^{\xi_1}) - \Delta_t(\mathcal{H}, \Psi^{\xi_2}))dt\leq B_2\|\xi_1 -\xi_2\|_{ham}.$$
\item Let $\Phi = \{\phi_t\}$ and $\Psi = \{\psi_t\}$ be two isotopies such that 
$\bar d(\Phi, \Psi)\leq r(g)/2$ where $r(g)$ is the injectivity radius of the Riemannian metric $g$ on $M$. Then, 
$$\int_0^1osc(\Delta_t(\mathcal{H},\Phi) - \Delta_t(\mathcal{H}, \Psi))dt\leq 4\max_t|\mathcal{H}_t|\bar d(\Phi, \Psi).$$
\end{enumerate}
\end{lemma}
{\it Proof.} For each $j = 1,2$, differentiating the reparameterized path
$\Psi^{\xi_j}$ in the variable $t$ yields 
$\dot\Psi^{\xi_j}(t) = \dot\xi_j(t)\dot\psi_{\xi_j(t)}$ for all $t\in[0,1].$ Compute   
$$\Delta_t(\mathcal{H},\Psi^{\xi_j}) = \int_0^{t}\mathcal{H}_t(\dot\Psi^{\xi_j}(s))\circ\Psi^{\xi_j}(s)ds = 
\int_0^{t}\dot\xi_j(s)\mathcal{H}_t(\dot\psi_{\xi_j(s)})\circ\psi_{\xi_j(s)}ds,$$
for each $t$. By a suitable change of variable, the right-hand side in the above estimates is written as  
$ \int_0^{t}\dot\xi_j(s)\mathcal{H}_t(\dot\psi_{\xi_j(s)})\circ\psi_{\xi_j(s)}ds =\int_0^{\xi_j(t)}\mathcal{H}_t(\dot\psi_u)\circ\psi_udu$. 
 This in turn yields 
$$\int_0^1osc(\Delta_t(\mathcal{H},\Psi^{\xi_1}) - \Delta_t(\mathcal{H}, \Psi^{\xi_2}))dt = 
\int_0^1osc(\int_{\min\{\xi_1(t), \xi_2(t)\}}^{\max\{\xi_1(t), \xi_2(t)\}}\mathcal{H}_t(\dot\psi_u)\circ\psi_udu)dt
$$ $$
\leq2\sup_{s,t,x}|\mathcal{H}_t(\dot\psi_s)(x)|\|\xi_1 -\xi_2\|_{C^0}
$$ 

$$\leq2\sup_{s,t,x}|\mathcal{H}_t(\dot\psi_s)(x)|\|\xi_1 -\xi_2\|_{ham}.$$
 Therefore, the desired $B_2$ is given by    
  $B_2 = 2\sup_{s,t,x}|\mathcal{H}_t(\dot\psi_s)(x)|\textless\infty.$\\
For (2), we set $\Phi = \{\phi_t\}$ and $\Psi = \{\psi_t\}$, and from the bi-invariance of the metric $\bar d$, we derive from the assumption 
that $\bar d(\Phi,\Psi) = \bar d(\Phi\circ\Psi^{-1}, Id)\leq r(g)/2.$
Under this condition, it follows from the lines of proof of Lemma 3.2 found in \cite{BanTchu} that 
$$\int_0^1osc(\Delta_t( \mathcal{H} , \Phi\circ\Psi^{-1}))dt
\leq 4\max_t|\mathcal{H}_t|\bar d(\Phi\circ\Psi^{-1}, Id).$$ On the other hand, we derive from Lemma \ref{cP11} that 
$$\Delta_t( \mathcal{H} , \Phi\circ\Psi^{-1}) = \Delta_t( \mathcal{H} ,\Psi^{-1})  + \Delta_t( \mathcal{H} ,\Phi)\circ\Psi^{-1} + cte_1,$$
and 
$$0 = \Delta_t( \mathcal{H} ,\Psi^{-1})  + \Delta_t( \mathcal{H} ,\Psi)\circ\Psi^{-1} + cte_2.$$
That is, for each $t$, there exists a constant $C$ which depends on $t$ such that 
$$\Delta_t( \mathcal{H} , \Phi\circ\Psi^{-1}) = -\Delta_t( \mathcal{H} ,\Psi)\circ\Psi^{-1}  + \Delta_t( \mathcal{H} ,\Phi)\circ\Psi^{-1} + C,$$
i.e. 

$$\int_0^1osc(-\Delta_t( \mathcal{H} ,\Psi)\circ\Psi^{-1}  + \Delta_t( \mathcal{H} ,\Phi)\circ\Psi^{-1} )dt 
= \int_0^1osc(\Delta_t( \mathcal{H} , \Phi\circ\Psi^{-1}))dt
$$ $$\leq 4\max_t|\mathcal{H}_t|\bar d(\Phi\circ\Psi^{-1}, Id).$$
 That achieves the proof.$ \hspace{0.1cm}\Box$\\

 {\it Proof of Lemma \ref{4ch2l4}.} Assume that 
$\Phi$ is generated by $( U,\mathcal{H})$, and consider $\Xi$ to be the path 
obtained by regularizing the path $\Phi$ as explained in Section \ref{SC7}. It follows 
that the path $\Xi$ is generated by an element $(V,\mathcal{K})$ 
so that 
$$l^{(1,\infty)}(\Xi) = \int_0^1 (osc(V_t) + |\mathcal{K}_t|)dt
\leq l^{(1,\infty)}(\Phi) + \int_0^1osc( r_t)dt + \int_0^1 osc(\Delta_t(\mathcal{H}, \phi_{(r,0)}))dt,$$
where $\phi_{( r,0)}$ is a Hamiltonian loop such that $\int_0^1osc( r_t)dt\textless \epsilon/2$ (see Section \ref{SC7} of the present paper). 
Since Polterovich's arguments provided in Section \ref{SC7} state that  
the path $\phi_{( r,0)}$ is arbitrarily close to the constant path identity (in the $C^\infty-$topology), we then derive from 
Lemma \ref{S1C42}-(2) that $\int_0^1osc(\Delta_t(\mathcal{H}, \phi_{( r,0)}))dt \textless \epsilon/2.$ We summarize the 
 above statements to get $ l^{(1,\infty)}(\Xi) \leq l^{(1,\infty)}(\Phi)  + \epsilon.$ 
Now, we use the path $\Xi$ to define a curve $\zeta$ as explained in Section \ref{SC7} of the present paper. Let $\Xi^{\zeta}$ be the path 
obtained by a reparameterization of $\Xi$ via $\zeta$. 
For each $s$, set $ \Omega_s =  \zeta'(s)(osc(V_{\zeta(s)}) + |\mathcal{K}_{\zeta(s)}|),$ 
 and derive from  (\ref{FNB}) that $ \Omega_s = l^{(1,\infty)}(\Xi)$. This 
yields  
 $\max_s\Omega_s = l^{\infty}(\Xi^\zeta) = l^{(1,\infty)}(\Xi)$. It follows from the above arguments that  
$l^{\infty}(\Xi^\zeta) = l^{(1,\infty)}(\Xi)$ and $ l^{(1,\infty)}(\Xi)  \textless l^{(1,\infty)}(\Phi) + \epsilon$, i.e.   
 $ l^{\infty}(\Xi^\zeta) \textless l^{(1,\infty)}(\Phi) + \epsilon.$ 
Therefore, to complete the proof, it suffices to take $\Psi = \Xi^\zeta$. 
$\Box$\\

The above proof of Lemma \ref{4ch2l4} is a refinement and a simplification of the proof of similar result found in \cite{BanTchu}. This 
is done in view to facilitate the readers to better understand the result of this paper. 

\begin{remark} As in \cite{Oh-M07}, 
  we have $l^{(1,\infty)}(.) \leq l^{\infty}(.) $ in general, where the former is invariant under reparameterization, 
while the latter is far from being invariant. But, as we can read in the proof of Lemma \ref{4ch2l4},  
any regular symplectic path $\Phi$ can be reparameterized to obtain another path $\Psi$ with the same extremities 
than $\Phi$ so that  $l^{(1,\infty)}(\Psi) = l^{\infty}(\Psi)$. 
\end{remark}
\subsection{Some remarks on flux geometry}

Another main result of this paper deals with flux geometry.\\

 It is known that the symplectic area swept by any smooth loop in $M$ 
under the action of any Hamiltonian isotopy is null:

 \begin{theorem}\label{Bn}(Banyaga, \cite{Ban78}).  
  Let $(M,\omega)$ be a closed symplectic manifold and $\Phi$ 
be any Hamiltonian isotopy. For any loop $\gamma\subset M$, we have 
$ Flux(\Phi).[\gamma] = 0.$ 
 \end{theorem}

As far as i known, 
an explicit study of the dual (in a certain sense) of the latter result is 
not yet done. This can be formalized as follows. 
 Given an arbitrary non-Hamiltonian isotopy, is there any property that may satisfy a non trivial loop in $M$ 
so that the symplectic area swept by the latter under the isotopy in question vanishes?   
More generally, given any non exact closed $1-$form $\alpha$ over $M$, is there any constructive method for generating non  
trivial smooth paths $\gamma$ in $M$ such that $\int_{\gamma}\alpha = 0$?\\

It is not too hard to see that for any given closed $1-$form 
$ \alpha,$ the function $\Delta_1(\alpha,\Phi)$ has a precise geometrical meaning: 
for each $x\in M$, where $\gamma_{x,\Phi}(t) = \Phi_t(x),$ i.e. for each $x\in M$, we can describe the real number 
 $\Delta_1(\alpha,\Phi)(x) $ as the algebraic value of the symplectic area of the $2-$chain swept by the orbit $t\longmapsto \Phi_t(x),$
 under the symplectic flow generated by $\alpha$. It follows 
quite naturally from the above arguments that each zero of the function $\Delta_1(\alpha,\Phi) $ 
gives rise to a null symplectic area $2-$chain or a solution of the equation $$\int_{\gamma}\alpha = 0,$$
with unknown 
$\gamma$. However, we have no guarantee whether such a function always admit at least a zero or not. Here is a sufficient condition 
which guarantees the existence of at least one zero for such a function.

\begin{lemma}\label{nullT} (Hamiltonian criterion). Let $(M,\omega)$ be a closed symplectic manifold. Let $\alpha$ 
be a closed $1-$form over $M$ , and  $\Phi$ be a Hamiltonian isotopy. Then for each representative $\Psi$
in the homotopic class of  $\Phi$ (relatively to fix extremities) 
the function $x\longmapsto \Delta_1(\alpha,\Psi)(x)$ has at least one zero in $M$.
\end{lemma}

{\sl Proof~:} Since the function $x\longmapsto \Delta_1(\alpha,\Phi)(x)$  is smooth and $M$ compact, the latter function achieves its 
bounds. This suggests that    
$$ \min_{x\in M}\Delta_1(\alpha,\Phi)(x)\int_M\omega^n\leq 
 \int_M \Delta_1(\alpha,\Phi)\omega^n\leq \max_{x\in M} \Delta_1(\alpha,\Phi)(x)\int_M\omega^n,$$ 
i.e.  
 $$\min_{x\in M}\Delta_1(\alpha,\Phi)(x)\leq 0,$$ and $$0\leq \max_{x\in M}\Delta_1(\alpha,\Phi)(x),$$  
since $\Phi$ is Hamiltonian (see Proposition \ref{SC41} of the present paper). 
On the other hand,
consider the following Poincar\'e's scalar product :
$$\langle,\rangle_P: H^1(M,\mathbb{R})\times H^{2n-1}(M,\mathbb{R})\rightarrow \mathbb{R}, $$ $$
([\alpha], [\beta])\mapsto \int_{M}\alpha\wedge\beta,$$
where $H^\ast(M,\mathbb{R}) $ represents the $\ast-$th de Rham cohomology group. 
Using this bilinear mapping, one checks from the proof of Proposition \ref{SC41} that 
\begin{equation}\label{EQT}
 \int_M \Delta_1(\alpha,\Phi)\omega^n = n\langle Flux(\Phi), [\alpha\wedge\omega^{(n-1)}]\rangle_P
\end{equation}
for all symplectic isotopy $\Phi,$ where $Flux$ represents the first Calabi's invariant. 
In particular, we see that the mean value $ \int_M \Delta_1(\alpha,\Phi)\omega^n$ depends only on the homotopic class 
of of $\Phi$ (relatively to fix extremities). Therefore, we get 
$$\int_M \Delta_1(\alpha,\Psi)\omega^n = \int_M \Delta_1(\alpha,\Phi)\omega^n = 0.$$ 

 This completes the proof. $\Box$\\ 

The third main result of this paper is the following theorem. 
 
\begin{theorem}\label{nullCH} Let $(M,\omega)$ be a closed symplectic manifold. Let $\Psi$ 
be any symplectic isotopy whose flux is not trivial. Then 
 any loop $\gamma\subset M$ which is homotopic to a closed Hamiltonian orbit (relatively to a fixed base point) trivializes the flux of $\Psi$, i.e.  $ Flux(\Psi).[\gamma] = 0.$
\end{theorem}
Theorem \ref{nullCH} states that on a closed symplectic manifold, if a loop is homotopic  to 
a closed Hamiltonian orbit (relatively to a fixed base point), then the symplectic area swept by the latter under the symplectic 
flow generated by any non exact closed $1-$form 
is trivial. This can be viewed in a certain sense as the dual form   
of Theorem \ref{Bn}.\\

{\sl Proof of Theorem \ref{nullCH}~:} Let $\Psi$ be a symplectic isotopy whose flux is not trivial, and denote by  
 $\mathcal{H}^\Psi$ the 
harmonic representative of the de Rham cohomology class  $ Flux(\Psi)$. Consider a 
Hamiltonian loop $\Phi = (\phi_t)$ in the fundamental group $\pi_1(Symp_0(M,\omega))$ of the group $Symp_0(M,\omega)$. In particular, 
since $\Phi$ is a loop we derive from Equation (\ref{Sholum}) that $$d\Delta_1(\mathcal{H}^\Psi,\Phi) = 0,$$ i.e. the function 
$\Delta_1(\mathcal{H}^\Psi,\Phi)$ is constant. But, Lemma \ref{nullT} suggests that the latter function must 
vanish since $\Phi$ is Hamiltonian. Therefore, the connectedness of $M$ imposes that the function $\Delta_1(\mathcal{H}^\Psi,\Phi)$ 
must be trivial. For each $x\in M$, consider the loop $\gamma_{x,\Phi}^t = \phi_t(x),$ and check that
$$ Flux(\Psi).[\gamma_{x,\Phi}] = \int_{\gamma_{x,\Phi}}\mathcal{H}^\Psi = 
\int_0^1\mathcal{H}^\Psi_{\gamma_{x,\Phi}^s}(\dot\gamma_{x,\Phi}^s)ds 
  = \Delta_1(\mathcal{H}^\Psi,\Phi)(x) = 0.$$ 
Now, let $\beta$ be a representative in $ [\gamma_{x,\Phi}],$ and let $h_x$ denotes the homotopy between $\beta$ and $\gamma_{x,\Phi}$. 
Consider the smooth $2-$chain $ \oplus(\beta, \gamma_{x,\Phi}) := \{h_x(s,t): 0\leq s,t\leq 1\}$. Since the harmonic $1-$form $\mathcal{H}^\Psi$ 
is closed, it follows from Stokes' theorem that 
$$0= \int_{ \oplus(\beta, \gamma_{x,\Phi})}d\mathcal{H}^\Psi = \int_{ \partial\oplus(\beta, \gamma_{x,\Phi})}\mathcal{H}^\Psi = \int_{\gamma_{x,\Phi}}\mathcal{H}^\Psi - \int_{\beta}\mathcal{H}^\Psi.$$
This completes the proof. $\Box$\\

Theorem \ref{nullCH} states that any loop $\gamma\subset M$ 
which is homotopic (relatively to a fixed base point) to a closed Hamiltonian orbit satisfies $\int_\gamma\alpha = 0,$  
for all closed $1-$form $\alpha$. This seems to suggest that each Hamiltonian loop can be viewed as the trivial element 
in $Hom(H^1(M,\mathbb{R}), \mathbb{R})$. 
More generally, this tells us that there is a linear and continuous mapping 
$$ \widetilde{\mathcal{K}}(\Phi) : H^1(M,\mathbb{R})\rightarrow\mathbb{R}, [\alpha]\longmapsto \dfrac{1}{n}\int_{M}\Delta_1(\alpha,\Phi)\omega^n,$$
for all fixed isotopy $\Phi$. That is,  $\widetilde{\mathcal{K}}(\Phi)$ belongs to $Hom(H^1(M,\mathbb{R}), \mathbb{R})$
 which is isomorphic to  $ H_1(M,\mathbb{R})$. So, there is a natural map:
$$\widetilde{\mathcal{K}}: Iso(M,\omega)\rightarrow H_1(M,\mathbb{R}), \Phi \longmapsto \widetilde{\mathcal{K}}(\Phi).$$
This symplectic invariant looks to be very similar to the usual mass flow.

\section{Some auxiliary results}
 In this section we prove the symplectic analogues of 
 Lemma 3.20, Lemma 3.21 and the $L^{(1,\infty)}-$ approximation found in \cite{Oh-M07}.

\begin{definition}(\cite{Ban08a})
 The $L^{(1,\infty)}$ Banyaga's topology on the space $Iso(M,\omega)$ is 
the metric topology induced by the following metric :
\begin{equation}\label{metri1}
  D^1 ((U,\mathcal{H}), (V ,\mathcal{K})) = \dfrac{D_0((U,\mathcal{H}), (V ,\mathcal{K})) + 
D_0(\overline{(U,\mathcal{H})}, \overline{(V ,\mathcal{K})})}{2},
\end{equation}
where
\begin{equation}\label{1metri1}
 D_0((U,\mathcal{H}), (V ,\mathcal{K})) = \int_0^1 osc(U_t - V_t) + |\mathcal{H}_t - \mathcal{K}_t|dt.
\end{equation}
 
\end{definition}
We will need the following lemma. 
\begin{lemma}\label{S1C41}
 Let $(M, g)$ be a closed  oriented Riemannian manifold. Let $\mathcal{H}\in  \mathcal P\mathcal{H}^1(M,g)$, 
let $\Phi= \{\phi_t\}$ be an isotopy, 
and let $\xi_1,\xi_2 :[0,1]\rightarrow[0,1] $ be two smooth functions such that $\xi_1$ is monotonic. Then 
 there exists a constant $B_1$ which depends on $\mathcal{H}$ and $\Phi$ such that 
 $$\int_0^1osc(\Delta_t(\mathcal{H}^{\xi_1},\Phi) - \Delta_t(\mathcal{H}^{\xi_2}, \Phi))dt\leq B_1\|\xi_1 -\xi_2\|_{ham}.$$

\end{lemma}
{\it Proof}. Since 
$\Delta(\mathcal{H}^{\xi_2}, \Phi) - \Delta(\mathcal{H}^{\xi_1}, \Phi) = 
\Delta_t(\mathcal{H}^{\xi_2} -\mathcal{H}^{\xi_1}, \Phi),$
we derive from (\ref{Mos})  that 

\begin{equation}\label{HRT}
 \int_0^1 osc(\Delta_t(\mathcal{H}^{\xi_2}, \Phi) - \Delta_t(\mathcal{H}^{\xi_1}, \Phi))dt 
\leq 2 diam(M)(1 + \sup_{t,s}|D\phi_{t}( \gamma_{y_0}(s))|) \int_0^1|\mathcal{H}^{\xi_1}_t - \mathcal{H}^{\xi_2}_t|dt,
\end{equation}
where $D\phi_t$ is the tangent map of $\phi_t$, $y_0$ a point in $ M$, and $ \gamma_{y_0}$ a  
minimizing geodesic such that $ \gamma_{y_0}(1) = y_0$ (see Section \ref{HOPFRINOW} of the present paper). Since  
$$ |\mathcal{H}^{\xi_1}_t - \mathcal{H}^{\xi_2}_t|  
\leq  |\dot\xi_1(t)\mathcal{H}_{\xi_1(t)}- \dot\xi_1(t)\mathcal{H}_{\xi_2(t)} | + 
|\dot\xi_1(t)\mathcal{H}_{\xi_2(t)}- \dot\xi_2(t)\mathcal{H}_{\xi_2(t)}|,$$ we use the Lipschitz nature of 
 the map $t\mapsto \mathcal{H}_t$ to derive the existence of a constant $c_0\textgreater 0$ which depends on $\mathcal{H}$ such that, 
$$ |\mathcal{H}^{\xi_1}_t - \mathcal{H}^{\xi_2}_t| \leq \max_{t}|\mathcal{H}_t||\dot\xi_1(t) - \dot\xi_2(t)|
+ c_0\|\xi_1 - \xi_2\|_{C^0}|\dot\xi_1(t)|,$$ for each $t$, and this yields,
$$\int_0^1|\mathcal{H}^{\xi_1}_t - \mathcal{H}^{\xi_2}_t|dt\leq\max_{t}|\mathcal{H}_t|\int_0^1|\dot\xi_1(t) - \dot\xi_2(t)|dt
+ c_0\|\xi_1 - \xi_2\|_{C^0}$$ $$\leq2 \max(c_0 ,\max_{t}|\mathcal{H}_t|)\|\xi_1 - \xi_2\|_{ham}.$$
Inserting the above estimates  in  (\ref{HRT}) yields
 
$$  \int_0^1 osc(\Delta_t(\mathcal{H}^{\xi_2}, \Phi) - \Delta_t(\mathcal{H}^{\xi_1}, \Phi))dt\leq B_1\|\xi_1 - \xi_2\|_{ham},$$
where    
$B_1 = 4 diam(M)\max(c_0 ,\max_{t}|\mathcal{H}_t|)(1 + \sup_{t,s}| D\phi_{t}(\gamma_{y_0}(s))|)\textless\infty.$
This completes the proof.$ \hspace{8.8cm}\Box$

\begin{lemma}\label{RL2}
 Let $(M,\omega)$ be a closed  
symplectic manifold. Let $( U,\mathcal{H})\in\mathfrak{T}(M, \omega, g)$, and $\xi_1, \xi_2 : [0,1] \rightarrow[0,1]$ 
be two smooth functions. Assume that $\xi_1$ is monotone.
Then there exists a constant $C$ which depends on $(U,\mathcal{H})$ such that, 
$$D^1(( U,\mathcal{H})^{\xi_1}, (U,\mathcal{H})^{\xi_2})\leq C\|\xi_1 -\xi_2\|_{ham}.$$ 
\end{lemma}
We shall give a complete proof of this lemma later on. The following result is an immediate consequence of Lemma \ref{RL2}.
\begin{lemma}\label{RL3}
 Let $(M,\omega)$ be a closed  
symplectic manifold. Let $(U^i,\mathcal{H}^i)$ be a Cauchy sequence in $D^1$, 
and $\xi_1, \xi_2 : [0,1] \rightarrow[0,1]$ be two monotone smooth functions. Given $\epsilon\textgreater0$, there exists two positive constants 
$\delta = \delta(\{(U^i,\mathcal{H}^i)\})$, and $j_0 = j_0(\{(U^i,\mathcal{H}^i)\}),$ such that : if $\xi_1, \xi_2$ 
satisfy $\|\xi_1 -\xi_2\|_{ham}\textless\delta$, 
then $$D^1(( U^i,\mathcal{H}^i)^{\xi_1}, (U^i,\mathcal{H}^i)^{\xi_2})\textless\epsilon,$$ for all $i\geq j_0$.
\end{lemma}

{\it Proof.} 
Since the sequence $(U^i,\mathcal{H}^i)$ is Cauchy in $D^1$,  
one can find an integer $j_0$ large enough such that $D^1(( U^i,\mathcal{H}^i)^{\xi_1}, 
(U^{j_0},\mathcal{H}^{j_0})^{\xi_1})\textless\epsilon/3$ for all $i\geq j_0$. 
Next, we choose $\delta = \epsilon/3C$ where $C$  
is obtained by applying Lemma \ref{RL2} to $(U^{j_0},\mathcal{H}^{j_0})$, and derive that
$$D^1(( U^i,\mathcal{H}^i)^{\xi_1}, (U^i,\mathcal{H}^i)^{\xi_2})\leq D^1(( U^i,\mathcal{H}^i)^{\xi_1}, (U^{j_0},\mathcal{H}^{j_0})^{\xi_1})$$
 $$ + D^1(( U^{j_0},\mathcal{H}^{j_0})^{\xi_1}, (U^{j_0},\mathcal{H}^{j_0})^{\xi_2})+
D^1(( U^{j_0},\mathcal{H}^{j_0})^{\xi_2}, (U^{i},\mathcal{H}^{i})^{\xi_2})$$ $$\leq \epsilon/3 + \epsilon/3 + \epsilon/3,$$
as long as $\|\xi_1 -\xi_2\|_{ham}\textless\delta$, and $i\geq j_0.$ This completes the proof.$ \hspace{2.1cm}\Box$\\

The following result is the symplectic analogue of a slight variation of the $L^{(1,\infty)}-$ approximation lemma found in \cite{Oh-M07}. 
It shows that any symplectic isotopy can be approximated arbitrarily closely in $ D^1$ by a boundary flat 
symplectic path (with the same extremities) so that they are also close 
in $\bar d$.
\begin{lemma}\label{L2}
 Let $(M,\omega)$ be a closed  
symplectic manifold. Let $\Phi = \phi_{( U,\mathcal{H})}$ be a symplectic isotopy, and let $\epsilon$ be a positive real number. Then, 
there exists a boundary flat symplectic isotopy $\Psi = \psi_{( V,\mathcal{K})}$ with the same extremities than $\Phi$ such that 
$D^1(( U,\mathcal{H}), ( V,\mathcal{K}))\textless\epsilon,$ and  $\bar{d}(\Psi, \Phi) \textless \epsilon.$
\end{lemma}

{\it Proof.} Let $\epsilon$ be a positive real number. We consider $\xi :[0,1]\rightarrow[0,1] $ to be any smooth positive and 
increasing function, which is constant on the intervals $[0,  \delta]$ and $[1 - \delta, 1]$ where 
 $0\textless\delta\textless 1/13$. Next, we define $ (V,\mathcal{K})$ to be the element 
$(U,\mathcal{H})^{\xi}$ as explained in Section \ref{RPSC4}.  
It follows from the definition of the curve   
$\xi$ that the symplectic isotopy $ \psi_{( V,\mathcal{K})}$ is boundary flat and has the same extremities than $\phi_{(U,\mathcal{H})} $. 
Applying Lemma \ref{RL2} with $\xi_1 = id$ and $\xi_2 = \xi$, we deduce that 
$D^1(( U,\mathcal{H}),( V,\mathcal{K})) \leq C\|\xi - id\|_{ham} $ where $C$ is the constant in Lemma \ref{RL2} which 
depends only on $(U,\mathcal{H})$. On the other hand, since the maps 
$(t,x)\mapsto\phi_{(U,\mathcal{H})}^{t}(x)$  and $(t,x)\mapsto\phi_{(U,\mathcal{H})}^{-t}(x)$  are Lipschitz continuous, it turns out that 
there exists a constant $l_0\textgreater0$ which 
depends only on $(U,\mathcal{H})$ such that 
$\bar d(\phi_{(U,\mathcal{H})}, \psi_{(V,\mathcal{K})})\leq l_0\|\xi - id\|_{C^0}\textless l_0\|\xi - id\|_{ham}$. 
To complete the proof, it suffices to choose the curve $\xi$ so that $\|\xi - id\|_{ham}
\leq\min\{\epsilon/C; \epsilon/l_0;\epsilon\}$. This completes the proof.$ \hspace{11cm}\Box$\\

{\it Proof of Lemma \ref{RL2}.} In the following $\Phi$ represents the symplectic isotopy generated by $( U,\mathcal{H})$.  
\begin{itemize}
 \item Step (1). Consider the normalized function $V = U^{\xi_1} - U^{\xi_2},$ and compute 
$$|V_t| = |\dot\xi_1(t)U_{\xi_1(t)} -\dot\xi_2(t)U_{\xi_2(t)} |
\leq |\dot\xi_1(t)||U_{\xi_1(t)} - U_{\xi_2(t)}| + |U_{\xi_2(t)}||\dot\xi_1(t) - \dot\xi_2(t)|,$$ 
for each $t$. Since $M$ is compact, 
we use the Lipschitz nature of the smooth map $(t,x)\mapsto U_t(x)$ to derive the existence of  
a constant $k_0\textgreater0$ depending on $U$ such that $\max_{x\in M}|U_{t}(x) - U_{s}(x)| \leq k_0|t - s|$ for all $t,s\in [0,1]$. 
This yields $$ 0\leq\max_{x\in M}V_t(x)\leq k_0|\dot\xi_1(t)||\xi_1(t) - \xi_2(t)| + \max_x(U_{t}(x))|\dot\xi_1(t) - \dot\xi_2(t)|.$$
Similarly, we derive that 
$$0\leq-\min_{x\in M}V_t(x)\leq k_0|\dot\xi_1(t)||\xi_1(t) - \xi_2(t)| -\min_x(U_{t}(x))|\dot\xi_1(t) - \dot\xi_2(t)|.$$
 It follows straight from the above estimates that 
\begin{equation}\label{a12}
 \int_0^1osc(V_t)dt\leq 2k_0\max_t|\xi_1(t) - \xi_2(t)| + \max_t(osc(U_{t}))\int_0^1|\dot\xi_1(t) - \dot\xi_2(t)|dt.
\end{equation}
\item Step (2). We set $ \mathcal{K}= \mathcal{H}^{\xi_1} - \mathcal{H}^{\xi_2}$, and compute 
$$ |\mathcal{K}_t| = |\dot\xi_1(t)\mathcal{K}_{\xi_1(t)} - \dot\xi_2(t)\mathcal{K}_{\xi_2(t)}| \leq |\xi_1(t) - \xi_2(t)||\dot\xi_1(t)|
 + | \mathcal{H}_{\xi_2(t)}||\dot\xi_1(t) - \dot\xi_2(t)|,$$
for each $t$. The Lipschitz nature of the smooth map $t\mapsto \mathcal{H}_t$, tells us 
 that there exists a constant $c_0\textgreater0$ which depends on $\mathcal{H}$ such that $|\mathcal{H}_t - \mathcal{H}_s|\leq c_0|t - s|$ 
for all $s,t\in [0,1]$. This yields
\begin{equation}\label{a12E}
 |\mathcal{K}_t|\leq c_0|\xi_1(t) - \xi_2(t)||\dot\xi_1(t)| + | \mathcal{H}_{\xi_2(t)}||\dot\xi_1(t) - \dot\xi_2(t)|. 
\end{equation}
Integrating (\ref{a12E}) in the variable $t$ yields,
\begin{equation}\label{a112}
 \int_0^1| \mathcal{K}_t|dt 
\leq 2c_0\max_t|\xi_1(t) - \xi_2(t)| + \max_t| \mathcal{H}_t|\int_0^1|\dot\xi_1(t) - \dot\xi_2(t)|dt
\end{equation}
 Adding  (\ref{a12}) and  (\ref{a112}) together we get
\begin{equation}\label{a16}
 D_0(( U,\mathcal{H})^{\xi_1}, (U,\mathcal{H})^{\xi_2})\leq 4\max\{k_0 + c_0, \max_t(| \mathcal{H}_t| + osc(U_t)) \}\|\xi_1 -\xi_2\|_{ham}
\end{equation}
\item Step (3). On the other hand,  we compute  
$$\overline{( U,\mathcal{H})^{\xi_j}} = (-U^{\xi_j}\circ\Phi^{\xi_j} - \widetilde 
\Delta(\mathcal{H}^{\xi_j},\Phi^{\xi_j}),- \mathcal{H}^{\xi_j}),$$ 
for each $j = 1, 2$ 
and derive from (\ref{1metri1}) that 
$$D_0(\overline{( U,\mathcal{H})^{\xi_1}}, \overline{( U,\mathcal{H})^{\xi_2}})
\leq \int_0^1osc(\widetilde\Delta_t(\mathcal{H}^{\xi_1},\Phi^{\xi_1}) - 
\widetilde\Delta_t(\mathcal{H}^{\xi_2},\Phi^{\xi_2}))dt$$ $$+ 
  \int_0^1 osc(U^{\xi_2} _t  
- U^{\xi_1}_t) + |\mathcal{H}^{\xi_1}_t- \mathcal{H}^{\xi_2}_t|dt + \int_0^1osc(U^{\xi_1}_t\circ\Phi^{\xi_2}(t)  
- U^{\xi_1}_t\circ\Phi^{\xi_1}(t))dt,$$
$$\leq  \int_0^1osc(\widetilde\Delta_t(\mathcal{H}^{\xi_1},\Phi^{\xi_1}) - 
\widetilde\Delta_t(\mathcal{H}^{\xi_2},\Phi^{\xi_2}))dt + D_0(( U,\mathcal{H})^{\xi_1},( U,\mathcal{H})^{\xi_2})$$
$$ + k_1\|\xi_1 -\xi_2\|_{ham}.$$
In the above estimates, to obtain the quantity $k_1\|\xi_1 -\xi_2\|_{ham}$ we use the Lipschitz nature of the 
 maps $(x,t)\mapsto U_t(x)$, $(x,t)\mapsto\Phi^{-1}(t)(x)$, and $(x,t)\mapsto\Phi(t)(x)$
to derive the existence of a constant $k_1\textgreater 0$ depending on $\Phi$ such that 
$\int_0^1osc(\dot\xi_1(t)U_{\xi_1(t)}\circ\Phi^{\xi_1}(t)  
- \dot\xi_1(t)U_{\xi_1(t)}\circ\Phi^{\xi_2}(t))dt\leq k_1\|\xi_1 -\xi_2\|_{ham}.$ We derive from  
triangle inequality that  
$$\int_0^1osc(\widetilde\Delta_t(\mathcal{H}^{\xi_1},\Phi^{\xi_1}) - 
\widetilde\Delta_t(\mathcal{H}^{\xi_2},\Phi^{\xi_2}))dt \leq \int_0^1osc(\widetilde\Delta_t(\mathcal{H}^{\xi_1},\Phi^{\xi_1}) - 
\widetilde\Delta_t(\mathcal{H}^{\xi_1},\Phi^{\xi_2}))dt $$ $$+ 
  \int_0^1osc(\widetilde\Delta_t(\mathcal{H}^{\xi_1},\Phi^{\xi_2}) - 
\widetilde\Delta_t(\mathcal{H}^{\xi_2},\Phi^{\xi_2}))dt.$$ But, applying Lemma \ref{S1C41} to $\mathcal{H}$ and $\Phi^{\xi_2}$
yields $$\int_0^1osc(\widetilde\Delta_t(\mathcal{H}^{\xi_1},\Phi^{\xi_2}) - 
\widetilde\Delta_t(\mathcal{H}^{\xi_2},\Phi^{\xi_2}))dt
\leq B_1\|\xi_1 -\xi_2\|_{ham},$$ where the positive constant $B_1$ depends on 
 $\Phi$, while it follows from the proof of Lemma \ref{S1C42}-(1) that 
$$\int_0^1osc(\widetilde\Delta_t(\mathcal{H}^{\xi_1},\Phi^{\xi_1}) - 
\widetilde\Delta_t(\mathcal{H}^{\xi_1},\Phi^{\xi_2}))dt $$ $$= 
\int_0^1osc(\int_{\min\{\xi_1(t), \xi_2(t)\}}^{\max\{\xi_1(t), \xi_2(t)\}}\dot\xi_1(t)\mathcal{H}_{\xi_1(t)}(\dot\psi_u)\circ\psi_udu)dt
$$ $$
\leq2\sup_{s,t,x}|\mathcal{H}_t(\dot\psi_s)(x)|\|\xi_1 -\xi_2\|_{C^0}\int_0^1\dot\xi_1(t) dt
$$ $$\leq B_2\|\xi_1 -\xi_2\|_{ham},$$

 where  $B_2 = 2\sup_{s,t,x}|\mathcal{H}_t(\dot\psi_s)(x)|$ depends on
 $\Phi$. Hence, we derive from the above statements that 
$$ D_0(\overline{( U,\mathcal{H})^{\xi_1}}, \overline{( U,\mathcal{H})^{\xi_2}}) \leq
D_0(( U,\mathcal{H})^{\xi_1}, ( U,\mathcal{H})^{\xi_2}) + k_1\|\xi_1 -\xi_2\|_{ham}$$ 
$$ + B_1\|\xi_1 -\xi_2\|_{ham} + B_2\|\xi_1 -\xi_2\|_{ham}.$$
\item Step (4). Since  
$$D^1(( U,\mathcal{H})^{\xi_1},( U,\mathcal{H})^{\xi_2}) = 
D_0(\overline{( U,\mathcal{H})^{\xi_1}}, \overline{( U,\mathcal{H})^{\xi_2}}) + D_0(( U,\mathcal{H})^{\xi_1},( U,\mathcal{H})^{\xi_2})/2,$$
 we derive from step (2) and step (3) that
$$D^1(( U,\mathcal{H})^{\xi_1},( U,\mathcal{H})^{\xi_2}) \leq C\|\xi_1 -\xi_2\|_{ham},$$
where $C = B_1 + B_2 + k_1 + 4\max\{k_0 + c_0, \max_t(| \mathcal{H}_t| + osc(U_t)) \} \textless \infty.$ This completes 
the proof.$ \hspace{8cm}\Box$

\end{itemize}

\section{Applications}
From the topological point of view, Lemma \ref{S1C42}-(2)  
suggests that on a closed symplectic manifold any smooth family   
 of harmonic $1-$forms $\mathcal{H}$ gives rise to a nontrivial Hamiltonian path which is small in 
the Hofer norms. As a matter of fact, let $\epsilon$ be an arbitrary positive real number.
We define a $C^0-$neighborhood $\mathcal{W}(\epsilon,\mathcal{H}, r(g), Id)$ of the identity map by 
 $$\mathcal{W}(\epsilon,\mathcal{H}, r(g), Id): 
= \{\Psi\in Iso(M, \omega)|  \bar d(\Psi,Id)\leq\min(r(g), \epsilon/[4\max_t|\mathcal{H}_t| + 1])\}.$$
According to Polterovich ( Proposition $5.2.A$, \cite{Polt93}) 
the set $\mathcal{W}(\epsilon,\mathcal{H}, r(g), Id)$ contains at least a nontrivial 
Hamiltonian loop. On the other hand, Lemma \ref{S1C42}-(2) suggests that for all $\Psi\in\mathcal{W}(\epsilon,\mathcal{H}, r(g), Id)$ 
we have $ \int_0^1 osc(\widetilde\Delta_t(\mathcal{H},\Psi))dt\textless \epsilon.$ But, it follows from the 
proof of  Lemma \ref{S1C42}-(2) that 
we also have $\max_t osc(\widetilde\Delta_t(\mathcal{H},\Psi))\textless \epsilon.$ Hence,   
the Hofer norms of the 
Hamiltonian path generated by $\widetilde\Delta(\mathcal{H},\Psi)$ are bounded from above by $\epsilon$. Furthermore, the mapping 
$\Delta_{\mathcal{H}} :
\Psi \mapsto\widetilde\Delta(\mathcal{H},\Psi)$
maps continuously $\mathcal{W}(\epsilon,\mathcal{H}, r(g), Id)$
 into a $C^0-$neighborhood of the trivial function. This agrees with Lemma 3.2 found in \cite{BanTchu}.\\

The following result is a consequence of Theorem \ref{nullCH}. It compares symplectic paths to mechanical motions. In particular, it    
tells us how a closed Hamiltonian orbit behaves (or winds) in  $\mathbb{T}^{2}$.
\begin{lemma}\label{Permutorus} Consider the $2-$dimensional revolution 
torus $\mathbb{T}^2 $ equipped with its standard symplectic form $\omega$. Let $\alpha$ be a non exact closed $1-$form 
over $ \mathbb{T}^2 $. 
Consider $R$ to be the mechanical motion represented by a complete 
rotation about the principal axis of $\mathbb{T}^2$. Then, either $R$ cannot be represented by the symplectic flow generated by $\alpha$ 
or there
 is no meridian circle in $\mathbb{T}^2$  which is an orbit of a Hamiltonian loop over $(\mathbb{T}^2, \omega)$. 
\end{lemma}

{\sl Proof~:} Assume that $R$ can be represented by the symplectic flow $(\theta^t)$ generated by $\alpha$
and that there exists a
 Hamiltonian loop $\Phi = (\phi^t)$ whose an orbit is a meridian circle $C_0$ in $\mathbb{T}^2.$
 Then by assumption there exists a point $z\in M$ such that the path $C_0 : t\mapsto\phi^t(z)$ is a meridian 
circle in  $\mathbb{T}^2$. According 
to Lemma \ref{nullCH}, we must have, $\Delta_1(\alpha,\Phi)(z)= 0,$ since $\Phi$ is Hamiltonian. On the other hand, 

$$\Delta_1(\alpha,\Phi)(z) = \int_{C_0}(\int_0^1\iota(\dot\theta^t)\omega dt) 
= \int_0^1\int_0^1\omega(\dot\theta^t(C_0(s)), \dot C_0(s)) dtds 
= \int_{[0,1]\times[0,1]}(\Theta_{C_0})^\ast\omega, $$
where $$\Theta_{C_0} : [0,1]\times [0,1]\rightarrow M, (t,s)\mapsto \theta^t(\phi^s(z)).$$ Then, we see that 
 $\Delta_1(\alpha,\Phi)(z)$ is the algebraic value of the volume of the set\\  
$ \{\theta^s(\phi_t(z))| 0\leq t,s\leq 1\},$ which is nothing than $\mathbb{T}^2.$ That is, 
 $$0 =  \Delta_1(\alpha,\Phi)(z)= Vol(\mathbb{T}^2) \neq 0.$$
This is a contradiction. The claim follows. 
$\Box$

\begin{center}
 {\bf Acknowledgments:}
 \end{center}
\begin{center}
 Thanks to Hodge's theory for enabling us to further understand  Banyaga's topologies, and some implications of Polterovich's works 
in the study of symplectic dynamics.\\
\end{center}
\begin{center} 
 I would like to thank the referees for carefully reading an earlier draft of this paper and suggesting some constructive hints. 
\end{center}

\bibliographystyle{plain}

\begin{thebibliography}{99}
\bibitem{Ban} 
A. Banyaga, : {\em Sur la structure de diff\'{e}omorphismes qui 
pr\'{e}servent une forme symplectique,}
 comment. Math. Helv. {\bf 53} $(1978)$ pp $174-2227$.
 \bibitem{BanTchu}
 A. Banyaga and S. Tchuiaga : \emph{The group of strong symplectic homeomorphisms in $L^\infty$-metric,} 
Adv. Geom. {\bf14} (2014), no. 3, $523–539.$
\bibitem{Ban08a}
A. Banyaga : \emph{A Hofer-like metric on the group of symplectic diffeomorphisms,} 
Symplectic topology and measure preserving dynamical systems, $(2010)$, pp. $1–23.$
\bibitem{Ban10c}
A. Banyaga : \emph{ On the group of strong symplectic homeomorphisms,}
 C. R. Math. Acad. Sci. Paris {\bf346} $(2008)$, no. $15-16,$ $867–872.$ 
\bibitem{Ban78}
 A. Banyaga, : \emph{Sur la structure de diff\'{e}omorphismes qui 
pr\'{e}servent une forme symplectique,}
 comment. Math. Helv. {\bf 53} $(1978)$ pp $174-2227$.
\bibitem{Hirs76}
 M. Hirsch : \emph{Differential Topology,} Graduate Texts in Mathematics, no. 33, Springer
  Verlag, New York-Heidelberg. {\bf 3} $(1976)$ corrected reprint $(1994)$.
 \bibitem{Hofer90}
 H. Hofer : \emph{On the topological properties of symplectic maps,}
 Proc. Royal Soc. Edinburgh {\bf 115A}(1990), pp $25-38$.
\bibitem{Hof-Zen94}
 H. Hofer and E. Zehnder : \emph{Symplectic invariants and Hamiltonian dynamics,}
 Birkhauser Advanced Texts, Birkhauser Verlag $(1994)$..
\bibitem{Lal-McD95} 
F. Lalonde and D. McDuff :\textit{The geometry of symplectic energy,}
Ann. of Math. {\bf 141}(1995), 711-727.

\bibitem{Oh-M07}
Y-G. Oh and S. M\"{u}ller : \emph{The group of Hamiltonian homeomorphisms and $C^0$-symplectic topology,}
J. Symp. Geometry {\bf 5}$(2007)$ $167-225$.

\bibitem{Polt93}
L. Polterovich : \emph{The Geometry of the Group of Symplectic Diffeomorphism,}
Lecture in Mathematics ETH Z\"{u}rich, Birkh\"{a}user Verlag, Basel-Boston $(2001)$.

\bibitem{Tchuia03}
S. Tchuiaga  : \emph{Some Structures of the Group of Strong Symplectic Homeomorphisms,} 
Global Journal of Advanced Research on Classical and Modern Geometry. Vol.2, Issue 1, pp.36-49, (2013)

 
\bibitem{FW71}
F. Warner : \emph{Foundation of differentiable manifolds and Lie groups,}
Graduate Texts in Mathematics, vol. $94$, Springer-Verlag, New York, $(1983)$.

 
  



 






 

 

 
%


 

\end{thebibliography}

%
%

\end{document}